\documentclass[10pt]{article}
\usepackage{amssymb}
\usepackage{amsmath}
\usepackage{epsfig}

\normalsize

\newtheorem{theorem}{Theorem}[section]
\newtheorem{lemma}[theorem]{Lemma}

\newtheorem{definition}[theorem]{Definition}

\newtheorem{algorithm}[theorem]{Algorithm}

\newtheorem{example}[theorem]{Example}

\newcommand{\boproof}{\noindent {\bf Proof. }}
\newcommand{\eoproof}{\hspace*{\fill} $\square$ \vspace{5pt}}
\newcommand{\red}{\sqsubseteq}

\newcommand{\R}{\mathbb R}
\newcommand{\Q}{\mathbb Q}
\newcommand{\Z}{\mathbb Z}

\newcommand{\Orthant}{\mathbb O}

\newcommand{\FourTiTwo}{{\tt 4ti2}}
\newcommand{\Kw}[1]{\underline{#1}}

\DeclareMathOperator{\newt}{Newt}
\DeclareMathOperator{\supp}{supp}
\DeclareMathOperator{\conv}{conv}

\DeclareMathOperator{\nF}{normalForm}
\DeclareMathOperator{\S-vectors}{S-vectors}

\begin{document}
\title{Computation of Atomic Fibers of $\Z$-Linear Maps}
\author{Raymond Hemmecke\\ Otto-von-Guericke-University Magdeburg}
\date{}
\maketitle

\begin{abstract}
For given matrix $A\in\Z^{d\times n}$, the set
$P^I_{A,b}=\{z:Az=b,z\in\Z^n_+\}$ describes the preimage or fiber
of $b\in\Z^d$ under the $\Z$-linear map
$f_A:\Z^n_+\rightarrow\Z^d$, $x\mapsto Ax$. The fiber $P^I_{A,b}$
is called atomic, if $P^I_{A,b}=P^I_{A,b_1}+P^I_{A,b_2}$ implies
$b=b_1$ or $b=b_2$. In this paper we present a novel algorithm to
compute such atomic fibers. An algorithmic solution to appearing
subproblems, application to integer programming, and computational
examples are included as well.
\end{abstract}


\section{Introduction}
Decomposition of rational polyhedra is at the heart of several
interesting applications. However, there are different definitions
of decomposability depending on the application. These definitions
mainly differ in the treatment of the (integer) points of the
polyhedron.

The simplest notion is that of linear decomposition of polyhedra.
Two polyhedra $P,Q\subseteq\R^n$ are called {\em homothetic} if
$P=\lambda Q+t$ for some $\lambda>0$ and $t\in\R^n$. Here, a
polyhedron $P$ is called {\em indecomposable}, if any
decomposition $P=Q_1+Q_2$ implies that both $Q_1$ and $Q_2$ are
homothetic to $P$. It can be shown that there are only finitely
many indecomposable rational polyhedra that are not homothetic to
each other. For further details on this type of decomposition we
refer the reader for example to
\cite{Gruenbaum:67,Henk+Koeppe+Weismantel:00,Kannan+Lovasz+Scarf:93,
McMullen:73,Meyer:74,Smith:87}.

Let us now come to a bit more restrictive decomposition. Here we
consider only polyhedra of the form $\{x\in\R^n:Ax\leq b\}$ for
given matrix $A\in\Z^{d\times n}$ and varying $b\in\Z^d$. To
emphasize that we only consider integer right-hand sides, we say
that a polyhedron $P$ is {\em integrally indecomposable}, if any
decomposition $P=Q_1+Q_2$ (into polyhedra with integer right-hand
sides) implies that both $Q_1$ and $Q_2$ are homothetic to $P$.
This decomposition is more restrictive than the linear
decomposition, since only such polyhedra $Q_1$ and $Q_2$ are
allowed that have an {\em integer} right-hand side. In
\cite{Henk+Koeppe+Weismantel:00}, the authors showed that
properties like TDI-ness of each member of a family of systems
$Ax\leq b$, $b\in\Z^d$, can be concluded from TDI-ness of the
integrally indecomposable systems.

Another important application of integral decomposition of
polyhedra is that of factorizing a multivariate polynomial. Here,
one considers only polyhedra of the form $\{x\in\R^n:Ax\leq b\}$
for given matrix $A\in\Z^{d\times n}$ and varying $b\in\Z^d$,
where each polyhedron is integer, that is, where each polyhedron
has only integer vertices. The reason for this restriction is the
simple observation that the so-called Newton polytope
$\newt(f):=\conv(\alpha:x^\alpha\in\supp(f))$ associated to a
polynomial $f$ is integer by definition. Moreover, the relation
$f=gh$ among three polynomials $f$, $g$, and $h$ implies
$\newt(f)=\newt(g)+\newt(h)$, a theorem due to Ostrowski. For more
details on this subject see for example \cite{Salem+Gao+Lauder:04}
and the references therein.

In this paper, we will consider another notion of decomposability.
In contrast to considering all points of a polyhedron, we restrict
out attention to its set of lattice points. More formally, we
consider family of sets of lattice points
\[
P^I_{A,b}:=\{z:Az=b,z\in\Z^n_+\},
\]
where the matrix $A\in\Z^{d\times n}$ is kept fix and the
right-hand side vector $b\in\Z^d$ varies. Note that $P^I_{A,b}$ is
exactly the preimage or {\em fiber} of $b$ under the linear map
$f_A:\Z^n_+\rightarrow\Z^d$, $x\mapsto Ax$. We call a fiber
$P^I_{A,b}$ {\em indecomposable} or {\em atomic}, if
$P^I_{A,b}=P^I_{A,b_1}+P^I_{A,b_2}$ implies $b=b_1$ or $b=b_2$.
Note that $P^I_{A,b}=P^I_{A,b_1}+P^I_{A,b_2}$ means that {\em
every lattice point} of $P^I_{A,b}$ is the sum of a {\em lattice
point} of $P^I_{A,b_1}$ and a {\em lattice point} of
$P^I_{A,b_2}$. This is indeed a very strong condition, but again
it can be shown that there are only finitely many (nonempty)
atomic fibers for a given matrix $A$ \cite{Maclagan:99}. Atomic
fibers were used in \cite{Adams+Hosten+Loustaunau+Miller:99} to
construct strong SAGBI-bases and receive more recent attention in
the computation of minimal vanishing sums of roots of unity
\cite{Steinberger:04}. Note that atomic fibers are not only
minimal (with respect to decomposability) within the given family,
but also generate every fiber $P^I_{A,b}$ in this family as a
Minkowski sum $P^I_{A,b}=\bigoplus_{i=1}^{k} \alpha_iP^I_{A,b_i}$,
$\alpha_i\in\Z_+$.

Below, we present a new algorithm to compute atomic fibers. In
fact, this algorithm can also compute -- what we call -- {\em
extended} atomic fibers. For this, we call
\[
Q^I_{A,b}:=\{z:Az=b,z\in\Z^n\}
\]
an extended (integer) fiber and we call it atomic, if
$(Q^I_{A,b}\cap\Orthant_j)=
(Q^I_{A,b_1}\cap\Orthant_j)+(Q^I_{A,b_2}\cap\Orthant_j)$ holds for
all the $2^n$ orthants $\Orthant_j$ of $\R^n$, then $b=b_1$ or
$b=b_2$. Also here it can be shown that there are only finitely
many (nonempty) extended atomic fibers for a given matrix. Also
this very strong notion of decomposability has an application: the
set ${\cal H}_\infty$ constructed in \cite{Hemmecke:SIP2} for use
in two-stage stochastic integer programming is in fact the set of
extended atomic fibers of the family of extended fibers
\[
\{x,y:x=b,Tx+Wy=0,x\in\Z^m,y\in\Z^n\}
\]
where $T$ and $W$ are kept fix and where $b$ varies.

One algorithm to compute atomic fibers via certain standard pairs
was presented in \cite{Adams+Hosten+Loustaunau+Miller:99}. Via
this algorithm, the authors computed the atomic fibers of the
twisted cubic by hand(!), see Example \ref{Example: Atomic Fibers
of Twisted Cubic}. Recently, the computation of atomic fibers
appears a sub-problem  or are being investigated in the design of
telecommunication networks with uncertain demand and capacity.

Before we present our algorithm, we include in Section
\ref{Section: Solving Integer Programs via Atomic Fibers} an
application of atomic fibers: the integer program
$\min\{c^\intercal z:z\in P^I_{A,b}\}$ can be solved via the
optimal solutions of the integer programs whose right-hand sides
are those of the atomic fibers.

\section{Atomic Fibers}

Let us now start our treatment with a formal definition of atomic
fibers.

\begin{definition}
For $u,v\in\R^n$ we say that $u\red v$ if $u^{(j)}v^{(j)}\geq 0$
and $|u^{(j)}|\leq |v^{(j)}|$ for all components $j=1,\ldots,n$,
that is, $u$ belongs to the same orthant as $v$ and its components
are not greater in absolute value than the corresponding
components of $v$.

For $S,T,U\subseteq\R^n$ we say that $S=T\oplus U$ if for every
$z\in S$ there are $z_1\in T$ and $z_2\in U$ with $z_1,z_2\red z$
and $z=z_1+z_2$.

For $A\in\Z^{d\times n}$ and $b\in\Z^d$ we call
$P^I_{A,b}:=\{z:Az=b,z\in\Z^n_+\}$ a {\bf fiber} and
$Q^I_{A,b}:=\{z:Az=b,z\in\Z^n\}$ an {\bf extended fiber} of $A$.

We call $P^I_{A,b}$ {\bf atomic}, if $P^I_{A,b}$ cannot be written
as $P^I_{A,b}=P^I_{A,b_1}\oplus P^I_{A,b_2}$ with two other fibers
$P^I_{A,b_1}$ and $P^I_{A,b_2}$.

Analogously, we call $Q^I_{A,b}$ {\bf atomic}, if $Q^I_{A,b}$
cannot be written as $Q^I_{A,b}=Q^I_{A,b_1}\oplus Q^I_{A,b_2}$
with two other extended fibers $Q^I_{A,b_1}$ and $Q^I_{A,b_2}$.

By $F(A)$ and by $E(A)$ denote the atomic and the extended atomic
fibers of $A$, respectively.
\end{definition}
Note that for every integer matrix $A$ there are at most finitely
many (extended) atomic fibers. This fact is based on the following
nice theorem.

\begin{theorem} [Maclagan, \cite{Maclagan:99}] \label{Theorem:
Maclagan} Let ${\cal I}$ be an infinite family of monomial ideals
in a polynomial ring $k[x_1,\ldots,x_n]$. Then there must exist
ideals $I,J\in {\cal I}$ with $I\subseteq J$.
\end{theorem}

To apply this theorem in our situation of atomic fibers, associate
to the fiber $Q^I_{A,b}$ the monomial ideal $I_{A,b}:=\langle
x^{(z^+,z^-)}:Az=b,z\in\Z^n\rangle\subseteq
\Q[x_1,\ldots,x_{2n}]$. $Q^I_{A,b}$ is an extended atomic fiber if
and only if $I_{A,b}$ is inclusion-maximal among all these ideals.

For the treatment below we like to point out here that the
relations $P^I_{A,b}=Q^I_{A,b}\cap\R^n_+$ and $P^I_{A,b_1}\oplus
P^I_{A,b_2}=P^I_{A,b_1}+P^I_{A,b_2}\subseteq P^I_{A,b_1+b_2}$ are
easy to verify. Another simple fact to notice is that
$Q^I_{A,b_1+b_2}=Q^I_{A,b_1}+Q^I_{A,b_2}$ whenever
$Q^I_{A,b_1}\neq\emptyset$ or $Q^I_{A,b_2}\neq\emptyset$.

\begin{example} \label{Example: Atomic Fibers of Twisted Cubic}
In \cite{Adams+Hosten+Loustaunau+Miller:99}, it was shown how
atomic fibers could be used to construct strong SAGBI bases for
monomial subalgebra over principal ideal domains. As an example,
they computed the atomic fibers of the matrix
$A=\left(\begin{smallmatrix} 3 & 2 & 1 & 0\\ 0 & 1 & 2 & 3\\
\end{smallmatrix}\right)$ by hand via an approach different from
the one we present below.

In the table below, we list the right-hand sides and all (finitely
many) elements in these $18$ atomic fibers. (Using \FourTiTwo\
\cite{4ti2}, in which the algorithm presented below that computes
(extended) atomic fibers has been implemented, we can not only
verify this result, but also quickly find the $51$ extended fibers
associated to this matrix.)
\[
\begin{array}{ll}
(0,3) & \{(0,0,0,1)\}\\
(1,2) & \{(0,0,1,0)\}\\
(2,1) & \{(0,1,0,0)\}\\
(3,0) & \{(1,0,0,0)\}\\
(2,4) & \{(0,1,0,1),(0,0,2,0)\}\\
(3,3) & \{(1,0,0,1),(0,1,1,0)\}\\
(4,2) & \{(0,2,0,0),(1,0,1,0)\}\\
(3,6) & \{(1,0,0,2),(0,1,1,1),(0,0,3,0)\}\\
(4,5) & \{(0,2,0,1),(0,1,2,0),(1,0,1,1)\}\\
(5,4) & \{(1,1,0,1),(0,2,1,0),(1,0,2,0)\}\\
(6,3) & \{(2,0,0,1),(1,1,1,0),(0,3,0,0)\}\\
(4,8) & \{(0,2,0,2),(1,0,1,2),(0,1,2,1),(0,0,4,0)\}\\
(6,6) & \{(2,0,0,2),(0,3,0,1),(1,1,1,1),(1,0,3,0),(0,2,2,0)\}\\
(8,4) & \{(2,1,0,1),(0,4,0,0),(1,2,1,0),(2,0,2,0)\}\\
(6,9) & \{(2,0,0,3),(0,3,0,2),(1,1,1,2),(1,0,3,1),(0,2,2,1),(0,1,4,0)\}\\
(9,6) & \{(3,0,0,2),(1,3,0,1),(2,1,1,1),(2,0,3,0),(1,2,2,0),(0,4,1,0)\}\\
(6,12) & \{(2,0,0,4),(0,3,0,3),(1,1,1,3),(1,0,3,2),(0,2,2,2),(0,1,4,1),(0,0,6,0)\}\\
(12,6) & \{(4,0,0,2),(2,3,0,1),(3,1,1,1),(3,0,3,0),(2,2,2,0),(0,6,0,0),(1,4,1,0)\}\\
\end{array}
\]
Thus, for example, the fiber given by the right-hand side $(8,7)$
is not atomic, since it can be decomposed into atomic fibers as
\[
P^I_{A,\left(\begin{smallmatrix}8\\7\\\end{smallmatrix}\right)}=
P^I_{A,\left(\begin{smallmatrix}2\\4\\\end{smallmatrix}\right)}\oplus
P^I_{A,\left(\begin{smallmatrix}6\\3\\\end{smallmatrix}\right)}.
\]
This can be quickly verified by looking at the elements in these
fibers:
\begin{eqnarray*}
& &
\{(2,1,0,2),(2,0,2,1),(1,1,3,0),(1,2,1,1),(0,4,0,1),(0,3,2,0)\}\\
& = &
\{(0,1,0,1),(0,0,2,0)\}\oplus\{(2,0,0,1),(1,1,1,0),(0,3,0,0)\}.
\end{eqnarray*}
Indeed, we have
\begin{eqnarray*}
(2,1,0,2) & = & (0,1,0,1)+(2,0,0,1), \\
(2,0,2,1) & = & (0,0,2,0)+(2,0,0,1), \\
(1,1,3,0) & = & (0,0,2,0)+(1,1,1,0), \\
(1,2,1,1) & = & (0,1,0,1)+(1,1,1,0), \\
(0,4,0,1) & = & (0,1,0,1)+(0,3,0,0), \\
(0,3,2,0) & = & (0,0,2,0)+(0,3,0,0).
\end{eqnarray*}
\end{example}
In Example \ref{Example: Atomic Fibers of Twisted Cubic} above, it
was easy to verify whether a given fiber is a summand in the
decomposition of another fiber by simply checking the finitely
many elements in the fiber for a decomposition. If the fibers are
not bounded, however, this would not give a finite procedure. The
following Lemma tells us how to solve this problem via the
(finitely many!) $\red$-minimal elements in the given fibers.

\begin{lemma}
Let $Q^I_{A,b_1}\neq\emptyset$ and $Q^I_{A,b_2}\neq\emptyset$.
Then $Q^I_{A,b_1+b_2}=Q^I_{A,b_1}\oplus Q^I_{A,b_2}$ if and only
if for every $\red$-minimal vector $v\in Q^I_{A,b_1+b_2}$ there is
a vector $w\in Q^I_{A,b_1}$ with $w\red v$.
\end{lemma}
\boproof Let $v\in Q^I_{A,b_1+b_2}$ and let $\bar{v}\in
Q^I_{A,b_1+b_2}$ be $\red$-minimal in $Q^I_{A,b_1+b_2}$ with
$\bar{v}\red v$. Thus, by the assumption in the lemma, there is
some $\bar{w}\in Q^I_{A,b_1}$ such that $\bar{w}\red\bar{v}$.
Clearly, this yields $\bar{v}-\bar{w}\in Q^I_{A,b_2}$ and
$\bar{v}-\bar{w}\red\bar{v}$.

We now claim that $v=(\bar{w}+v-\bar{v})+(\bar{v}-\bar{w})$ with
$\bar{w}+v-\bar{v}\in Q^I_{A,b_1}$, $\bar{v}-\bar{w}\in
Q^I_{A,b_2}$, $\bar{w}+v-\bar{v}\red v$, and $\bar{v}-\bar{w}\red
v$, is a desired representation of $v$. The first two relations
are trivial, if we keep in mind that $Av=A\bar{v}=b$,
$A\bar{w}=b_1$, and $b=b_1+b_2$. We get the other two relations as
follows:
\begin{itemize}
\item [(a)] $\bar{w}+v-\bar{v}\red \bar{v}+v-\bar{v}=v$, since by
construction $\bar{w}$ and $v-\bar{v}$ lie in the same orthant,
and
\item [(b)] $\bar{v}-\bar{w}\red\bar{v}\red v$, since
$\bar{w}\red\bar{v}$.
\end{itemize}

Thus, we have constructed for arbitrary $v\in Q^I_{A,b_1+b_2}$ a
valid representation of $v$ as a sum of two elements from
$Q^I_{A,b_1}$ and $Q^I_{A,b_2}$ that lie in the same orthant as
$v$. This concludes the proof. \eoproof

Using this lemma repeatedly, we are now able to find, for a given
right-hand side $b$, a decomposition
$Q^I_{A,b}=\bigoplus_{i=1}^{k}\alpha_iQ^I_{A,b_i}$,
$\alpha_i\in\Z_+$, that is, we can find a decomposition of a fiber
into a sum of atomic fibers. Note that replacing $Q$ by $P$ and
$E(A)$ by $F(A)$ in the algorithm below would yield a
decomposition of the atomic fiber: $P^I_{A,b}=\bigoplus_{i=1}^{k}
\alpha_iP^I_{A,b_i}$, $\alpha_i\in\Z_+$.

\begin{algorithm} \label{Algorithm to decompose fiber into sum of atomic fibers}
{(Algorithm to decompose extended fiber into extended atomic
fibers)}

\Kw{Input:} $A$, $E(A)=\{b_1,\ldots,b_k\}$

\Kw{Output:} $\alpha_1,\ldots,\alpha_k$ such that
$Q^I_{A,b}=\bigoplus\limits_{i=1}^{k} \alpha_iQ^I_{A,b_i}$

$\alpha_1:=\ldots:=\alpha_k:=0$

\Kw{for} i \Kw{from} 1 \Kw{to} k \Kw{do}

\hspace{1.0cm} \Kw{while} $Q^I_{A,b}=Q^I_{A,b_i}\oplus
Q^I_{A,b-b_i}$ \Kw{do}

\hspace{2.0cm} $b:=b-b_i$

\hspace{2.0cm} $a_i:=a_i+1$

\Kw{return} $\alpha_1,\ldots,\alpha_k$.
\end{algorithm}

It remains to state an algorithm that computes the finitely many
$\red$-minimal elements in $P^I_{A,b}$ and in $Q^I_{A,b}$,
respectively. One simple way to compute these elements is via
(truncated) Hilbert and Graver bases. For this we note that the
$\red$-minimal elements in $P^I_{A,b}$ correspond to those
elements in the {\em Hilbert basis} of the cone
$\{(z,u):Az-bu=0,(z,u)\in\R^{n+1}_+\}$ that have $u=1$, while the
$\red$-minimal elements in $Q^I_{A,b}$ correspond to those
elements $(z,u)$ in the {\em Graver basis} of the matrix $(A|-b)$
that have $u=1$. Finiteness of such Hilbert and Graver bases for
rational data $A,b$ immediately implies that there are indeed only
finitely many $\red$-minimal elements in $P^I_{A,b}$ and
$Q^I_{A,b}$, respectively. This, however, could have been already
concluded from the fact that $\red$-minimal elements in
$Q^I_{A,b}$ are in one-to-one correspondence with the (finitely
many) minimal generators of the monomial ideal $I_{A,b}$.

\section{Solving Integer Programs via Atomic Fibers}
\label{Section: Solving Integer Programs via Atomic Fibers}

Algorithm \ref{Algorithm to decompose fiber into sum of atomic
fibers} puts us in the position to solve integer programs with the
help of {\em atomic integer programs} whose right-hand sides
correspond to those of the atomic fibers. Therefore, these atomic
programs encode already the solutions to a whole family of integer
programs (for changing right-hand side).

\begin{definition} {(Atomic Programs)}
If $P^I_{A,b}$ is an atomic fiber of $A\in\Z^{d\times n}$, then we
call
\[
\min\{c^\intercal z:z\in P^I_{A,b}\}=\min\{c^\intercal
z:Az=b,z\in\Z^n_+\}
\]
the {\bf atomic} (integer) {\bf program} associated to
$P^I_{A,b}$.
\end{definition}

Once we know optimal solutions to these finitely many atomic
programs, we can solve the integer program for any given
right-hand side $b$ by simply decomposing the fiber $P^I_{A,b}$
into a sum of atomic fibers.

\begin{lemma} \label{Atomic fibers solve integer program}
Given $A\in\Z^{d\times n}$, $c\in\R^n$, and
$F(A)=\{b_1,\ldots,b_k\}$, let $z_1^{opt},\ldots,z_k^{opt}$ be
optimal solutions to
\[
\min\{c^\intercal z:Az=b_i,z\in\Z^n_+\},\;\;\; i=1,\ldots,k,
\]
and consider any solvable problem
\begin{equation}\label{Solvable Problem}
\min\{c^\intercal z:Az=b,z\in\Z^n\}.
\end{equation}
If $P^I_{A,b}=\bigoplus_{i=1}^{k} \alpha_iP^I_{A,b_i}$,
$\alpha_i\in\Z_+$, then $\bar{z}:=\sum_{i=1}^k\alpha_i z_i^{opt}$
is an optimal solution to (\ref{Solvable Problem}).
\end{lemma}

\boproof First we show that $\bar{z}$ is a feasible solution of
(\ref{Solvable Problem}). Clearly,
$\bar{z}=\sum_{i=1}^k\alpha_iz_i^{opt}\geq 0$. Moreover,
\[
A\bar{z}=A\left(\sum_{i=1}^k\alpha_iz_i^{opt}\right)
=\sum_{i=1}^k\alpha_i\left(Az_i^{opt}\right)
=\sum_{i=1}^k\alpha_ib_i=b.
\]
Therefore, $\bar{z}$ is a feasible solution of (\ref{Solvable
Problem}).

Let $z$ be any feasible solution of (\ref{Solvable Problem}). From
the decomposition $P^I_{A,b}=\bigoplus_{i=1}^{k}
\alpha_iP^I_{A,b_i}$ with $\alpha_i\in\Z_+$ we conclude that there
are $z_i\in P^I_{A,b_i}$, $i=1,\ldots,k$, such that
$z=\sum_{i=1}^k\alpha_iz_i$. But now we have
\[
c^\intercal z=c^\intercal\left(\sum_{i=1}^k\alpha_iz_i\right)=
\sum_{i=1}^k\alpha_i\left(c^\intercal z_i\right)\geq
\sum_{i=1}^k\alpha_i\left(c^\intercal z_i^{opt}\right)=
c^\intercal\left(\sum_{i=1}^k\alpha_i
z_i^{opt}\right)=c^\intercal\bar{z}.
\]
Therefore, $\bar{z}$ is optimal for (\ref{Solvable Problem}).
\eoproof

\section{Computation of Atomic Fibers}
\label{Section: Computation of Atomic Fibers}

In the following we show how to compute the finitely many
(extended) atomic fibers of a matrix. We use the algorithmic
pattern of a completion procedure. By $e_1,\ldots,e_n$ we denote
the unit vectors in $\R^n$.

\begin{algorithm} \label{Algorithm to compute atomic fibers}
{(Algorithm to compute atomic fibers)}

\Kw{Input:} $F:=\{Ae_1,\ldots,Ae_n,-Ae_1,\ldots,-Ae_n\}$

\Kw{Output:} a set $G$, such that $\{Q^I_{A,b}:b\in G\}$ contains
all extended atomic fibers of $A$

\vspace{0.2cm} $G:=F$

$C:=\bigcup\limits_{f,g\in G}\{f+g\}$ \hspace{4.2cm} (forming
$\S-vectors$)

\Kw{while} $C\neq \emptyset $ \Kw{do}

\hspace{1.0cm}$s:=$ an element in $C$

\hspace{1.0cm}$C:=C\setminus\{s\}$

\hspace{1.0cm}$f:=\nF(s,G)$

\hspace{1.0cm}\Kw{if} $f\neq 0$ \Kw{then}

\hspace{2.0cm}$G:=G\cup \{f\}$

\hspace{2.0cm}$C:=C\cup\bigcup\limits_{g\in G} \{f+g\}$
\hspace{2.2cm} (adding $\S-vectors$)

$G:=G\cup\{0\}$

\Kw{return} $G$.
\end{algorithm}

Behind the function $\nF(s,G)$ there is the following algorithm.

\begin{algorithm}
{(Normal form algorithm)}

\Kw{Input:} $s$, $G$

\Kw{Output:} a normal form of $s$ with respect to $G$

\vspace{0.2cm}

\Kw{while} there is some $g\in G$ such that
$Q^I_{A,s}=Q^I_{A,g}\oplus Q^I_{A,s-g}$ \Kw{do}

\hspace{1.0cm} $s:=s-g$

\Kw{return} $s$
\end{algorithm}

\begin{lemma} \label{Lemma: Algorithm to compute extended atomic fibers
terminates and is correct} Algorithm \ref{Algorithm to compute
atomic fibers} terminates and computes a set $G$ such that
$\{Q^I_{A,b}:b\in G\}$ contains all atomic fibers of $A$.
\end{lemma}

\boproof Associate with $b$ the monomial ideal $I_{A,b}:=\langle
x^{(z^+,z^-)}:Az=b,z\in\Z^n\rangle\subseteq
\Q[x_1,\ldots,x_{2n}]$. Algorithm \ref{Algorithm to compute atomic
fibers} generates a sequence $\{f_1,f_2,\ldots\}$ in $G\setminus
F$ such that $Q^I_{A,f_j}\neq Q^I_{A,f_i}\oplus Q^I_{A,f_j-f_i}$
whenever $i<j$. Thus, the corresponding sequence
$\{I_{A,f_1},I_{A,f_2},\ldots\}$ of monomial ideals satisfies
$I(A,f_j)\not\subseteq I(A,f_i)$ whenever $i<j$. We conclude, by
Maclagan's theorem \cite{Maclagan:99}, Theorem \ref{Theorem:
Maclagan} above, that this sequence of monomial ideals must be
finite and thus, Algorithm \ref{Algorithm to compute atomic
fibers} must terminate.

It remains to prove correctness. For this, let $G$ denote the set
that is returned by Algorithm \ref{Algorithm to compute atomic
fibers}. Moreover, let $Q^I_{A,\bar{b}}$ be an extended atomic
fiber of $A$ with $\bar{b}\neq 0$. We will show that $\bar{b}\in
G$.

Since $F\subseteq G\setminus\{0\}$, we know that
$Q^I_{A,\bar{b}}=\sum Q^I_{A,b_j}$ for finitely many (not
necessarily distinct) $b_j\in G\setminus\{0\}$. This implies in
particular, that every $z\in Q^I_{A,\bar{b}}$ can be written as a
sum $z=\sum v_j$ with $v_j\in Q^I_{A,b_j}$. We will show that we
can find $b_j\in G$ such that every $z\in Q^I_{A,\bar{b}}$ can be
written as a sum $z=\sum v_j$ with $v_j\in Q^I_{A,b_j}\subseteq
Q^I_{A,\bar{b}}$ and with $v_j\red z$. This implies
$Q^I_{A,\bar{b}}=\bigoplus Q^I_{A,b_j}$. Since $Q^I_{A,\bar{b}}$
is atomic, and thus indecomposable, this representation must be
trivial, that is, it has to be $Q^I_{A,\bar{b}}=Q^I_{A,\bar{b}}$,
and therefore we conclude $\bar{b}\in G$.

For $A\in\Z^{d\times n}$ and $b\in\Z^d$ let $M^I_{A,b}$ denote the
set of $\red$-minimal elements in $Q^I_{A,b}$. It is sufficient to
find $Q^I_{A,b_i}\in G$ such that every $z\in M^I_{A,\bar{b}}$ can
be written as a sum $z=\sum v_i$ with $v_i\in Q^I_{A,b_i}$, where
all $v_i$ lie in the same orthant as $z$ (as this then implies
$v_i\red z$). To see this, take any non-minimal $z'\in
Q^I_{A,\bar{b}}\setminus M^I_{A,\bar{b}}$ and let $z\in
M^I_{A,\bar{b}}$ be such that $z\red z'$. Then, clearly,
$h:=z'-z\in Q^I_{A,0}$ and every representation $z=\sum v_i$ with
$v_i\in Q^I_{A,b_i}$ and $v_i\red z$ readily extends to a desired
valid representation $z'=z+h=(v_{i_0}+h)+\sum_{i\neq i_0} v_i$,
since $v_{i_0}+h\in Q^I_{A,b_{i_0}}$ and since $v_{i_0}$ and $h$
and thus also $v_{i_0}+h$ lie in the same orthant as $z$ and $z'$.
Consequently, $h\red z'$, as required.

By the Gordan-Dickson Lemma, $M^I_{A,\bar{b}}$ is finite. Thus,
let denote its elements by $M^I_{A,\bar{b}}=\{z_1,\ldots,z_k\}$.
From all representations $Q^I_{A,\bar{b}}=\sum_{j\in J}
Q^I_{A,b_j}$ with $b_j\in G\setminus\{0\}$ and where
$z_i=\sum_{j\in J} v_{i,j}$ and $v_{i,j}\in Q^I_{A,b_j}$,
$i=1,\ldots,k$, choose one such that the sum
\begin{equation}\label{sum}
\sum_{i=1}^k\sum_{j\in J}\|v_{i,j}\|_1
\end{equation}
is minimal. By the triangle inequality we have that
\begin{equation}\label{triangle inequality}
\sum_{i=1}^k\sum_{j\in J}\|v_{i,j}\|_1\geq\sum_{i=1}^k \|z_i\|_1.
\end{equation}
Herein, equality holds if and only if all $v_{i,j}$ have the same
sign pattern as $z_i$, $i=1,\ldots,k$, that is, if and only if we
have $v_{i,j}\red z_i$ for all $i$ and all $j$. Thus, if we have
equality in (\ref{triangle inequality}) for such a minimal
representation $Q^I_{A,\bar{b}}=\sum_{j\in J} Q^I_{A,b_j}$, then
$v_{i,j}\in Q^I_{A,b_j}$ and $v_{i,j}\red z_i$ for all occurring
$v_{i,j}$, and we are done.

(It should be noted that we have required $b_j\in G\setminus\{0\}$
for all appearing $b_j$, that is in particular, $b_j\neq 0$. Those
$b_j$ will be sufficient to generate all $\red$-minimal elements
in the extended fiber $Q^I_{A,\bar{b}}$. The remaining elements in
$Q^I_{A,\bar{b}}$ we get by adding elements from $Q^I_{A,0}$ as
mentioned above.)

Therefore, let us assume that
\begin{equation}\label{inequality}
\sum_{i=1}^k\sum_{j\in J}\|v_{i,j}\|_1 > \sum_{i=1}^k z_i.
\end{equation}
In the following we construct a new representation
$Q^I_{A,\bar{b}}=\sum_{j'\in J'} Q^I_{A,b_j'}$ whose corresponding
sum (\ref{sum}) is smaller than the minimally chosen sum. This
contradiction proves that we have indeed equality in
(\ref{triangle inequality}) and our claim is proved.

From (\ref{inequality}) we conclude that there are indices
$i_0,j_1,j_2$ and a component $m\in\{1,\ldots,n\}$ such that
$v_{i_0,j_1}^{(m)}\cdot v_{i_0,j_2}^{(m)}<0$. The extended fiber
$Q^I_{A,b_{j_1}+b_{j_2}}$ was reduced to $Q^I_{A,0}$ by sets
$Q^I_{A,b_{j''}}$, $j''\in J''$. This gives representations
\[
v_{i,j_1}+v_{i,j_2}=\sum_{j''\in J''} w_{i,j''},\;\; w_{i,j''}\in
Q^I_{A,b_{j''}},w_{i,j''}\red v_{i,j_1}+v_{i,j_2}
\]
for $i=1,\ldots,k$. As all $w_{i,j''}$ lie in the same orthant of
$\R^n$ as $v_{i,j_1}+v_{i,j_2}$, we get
\[
\|\sum_{j''\in J''} w_{i,j''}\|_1 = \|v_{i,j_1}+v_{i,j_2}\|_1 \leq
\|v_{i,j_1}\|_1+\|v_{i,j_2}\|_1,
\]
with strict inequality for $i=i_0$.

Thus, replacing in $Q^I_{A,\bar{b}}=\sum_{j\in J} Q^I_{A,b_j}$ the
term $Q^I_{A,b_{j_1}}+Q^I_{A,b_{j_2}}$ by $\sum_{j''\in J''}
Q^I_{A,b_{j''}}$, we arrive at a new representation
$Q^I_{A,\bar{b}}=\sum_{j'\in J'} Q^I_{A,b_{j'}}$ whose
corresponding sum (\ref{sum}) is at most
\[
\sum_{i=1}^k\sum_{j'\in J'}\|v_{i,j'}\|_1 < \sum_{i=1}^k\sum_{j\in
J}\|v_{i,j}\|_1,
\]
contradicting the minimality of the representation
$Q^I_{A,\bar{b}}=\sum_{j\in J} Q^I_{A,b_j}$. This concludes the
proof. \eoproof

Having an algorithm available that computes all extended atomic
fibers, we can of course use it to compute atomic fibers: If
$P^I_{A,b}$ is atomic then so is $Q^I_{A,b}$, as any decomposition
of $Q^I_{A,b}$, restricted to the nonnegative orthant, would give
a decomposition of $P^I_{A,b}$. This way of computing atomic
fibers, however, is pretty slow, since there are far more extended
atomic fibers as atomic fibers. A similar behavior can be observed
when one extracts the Hilbert basis of the cone
$\{x:Ax=0,x\in\R^n_+\}$ from the Graver basis of $A$, as the
Graver basis is usually much bigger than the Hilbert basis one is
interested in.

In \cite{Hemmecke:Hilbert} it was shown that one can reduce this
difference in sizes by considering only a subset of $k$ variables,
by computing all indecomposable vectors in the set
$\{x:Ax=0,x\in\Z^k_+\times\Z^{n-k}\}$, and by iteratively adding a
new variable. The main observation exploited is that the
indecomposable vectors in $\{x:Ax=0,x\in\Z^k_+\times\Z^{n-k}\}$
are nonnegative on the first $k$ components and in addition
generate all indecomposable vectors in
$\{x:Ax=0,x\in\Z^{k+1}_+\times\Z^{n-k-1}\}$ as a nonnegative
integer linear combination. This keeps the sets from which the
indecomposable vectors are extracted comparatively small and leads
to a tremendous speed-up in the algorithm. With this algorithm,
bigger Hilbert bases, even with more than $500,000$ elements, can
be computed nowadays.

More interesting for our problem of finding atomic fibers is the
fact that a similar ``project-and-lift'' idea can be applied in
our situation here. Consider the set of ``partially extended''
fibers $Q^I_{A,b,k}:=\{z:Az=b,z\in\Z^k_+\times\Z^{n-k}\}$, where
$A$ is kept fix and where $b$ varies, and let us restrict the
conditions for decomposability to the first $k$ components. Now we
use Algorithm \ref{Algorithm to compute atomic fibers} to compute
a set of indecomposable/minimal partially extended fibers
$Q^I_{A,b_1,k},\ldots,Q^I_{A,b_s,k}$. They have the property of
being nonnegative on the first $k$ components in each vector and
of generating every other partially extended fiber $Q^I_{A,b,k}$
as a (nonnegative integer) Minkowski sum. This latter fact can be
used in the computation of indecomposable extended fibers
$Q^I_{A,b,k+1}$ by defining the input set to Algorithm
\ref{Algorithm to compute atomic fibers} to be
$F:=\{b_1,\ldots,b_s,Ae_{k-1},\ldots,Ae_n,-Ae_{k+1},\ldots,-Ae_n\}$.
Note that we do not need to add $-b_1,\ldots,-b_s$ here due to the
nice generating property that we have already achieved for
$Q^I_{A,b_1,k},\ldots,Q^I_{A,b_s,k}$.

Clearly, when $k=n$ we have $Q^I_{A,b,n}=P^I_{A,b}$ and we have
found the atomic fibers of $A$. Let us now demonstrate the
capability of this algorithm by considering the following example.
For this note that both Algorithm \ref{Algorithm to compute atomic
fibers} and the ``project-and-lift'' version of it to compute
atomic fibers directly are implemented in {\FourTiTwo}.

\begin{example} \label{Example: Atomic Fibers for 4x4}
One example that appears and was solved in \cite{Steinberger:04}
is the computation of the atomic fibers of the matrix
\[
\left(
\begin{array}{rrrrrrrrr}
1 & -1  &  0 & -1 &  1 &  0 &  0 &  0 & 0\\
0 &  1  & -1 &  0 & -1 &  1 &  0 &  0 & 0\\
0 &  0  &  0 &  1 & -1 &  0 & -1 &  1 & 0\\
0 &  0  &  0 &  0 &  1 & -1 &  0 & -1 & 1\\
\end{array}
\right).
\]
This matrix corresponds to a certain problem on $3\times 3$ tables
and has in fact $31$ atomic fibers and $79$ extended atomic
fibers, as can be easily verified with {\FourTiTwo}.

The next higher problem on $4\times 4$ tables leads to the matrix
\[
\left(
\begin{array}{rrrrrrrrrrrrrrrr}
1 & -1 &  0 &  0 & -1 &  1 &  0 &  0 &  0 &  0 &  0 &  0 &  0 &  0
& 0 & 0\\
0 &  1 & -1 &  0 &  0 & -1 &  1 &  0 &  0 &  0 &  0 &  0 &  0 &  0
& 0 & 0\\
0 &  0 &  1 & -1 &  0 &  0 & -1 &  1 &  0 &  0 &  0 &  0 &  0 &  0
& 0 & 0\\
0 &  0 &  0 &  0 &  1 & -1 &  0 &  0 & -1 &  1 &  0 &  0 &  0 &  0
& 0 & 0\\
0 &  0 &  0 &  0 &  0 &  1 & -1 &  0 &  0 & -1 &  1 &  0 &  0 &  0
& 0 & 0\\
0 &  0 &  0 &  0 &  0 &  0 &  1 & -1 &  0 &  0 & -1 &  1 &  0 &  0
& 0 & 0\\
0 &  0 &  0 &  0 &  0 &  0 &  0 &  0 &  1 & -1 &  0 &  0 & -1 &  1
& 0 & 0\\
0 &  0 &  0 &  0 &  0 &  0 &  0 &  0 &  0 &  1 & -1 &  0 &  0 & -1
& 1 & 0\\
0 &  0 &  0 &  0 &  0 &  0 &  0 &  0 &  0 &  0 &  1 & -1 &  0 &  0
& -1 & 1\\
\end{array}
\right).
\]
{\FourTiTwo} computes $12,675$ atomic fibers for this matrix.
\end{example}

\section{Conclusions}

In this paper we presented a completion procedure to compute
atomic and (extended) atomic fibers of the form
$P^I_{A,b}=\{z:Az=b,z\in\Z^n_+\}$ and
$Q^I_{A,b}=\{z:Az=b,z\in\Z^n\}$, respectively, and demonstrated
its ability to solve non-trivial examples.

If the number $n$ of variables is not fixed, the problem of
computing atomic fibers is NP-hard, since it is at least as hard
as computing Hilbert bases of rational cones, which are already
exponential in size. However, as Barvinok and Woods showed in
\cite{Barvinok+Woods:2003}, the Hilbert basis (and thus also the
Graver basis) can be encoded via a short (polynomial size)
rational generating function if the dimension is kept fixed. This
implies that the elements in each particular fiber
$P^I_{A,b}=\{z:Az=b,z\in\Z^n_+\}$ (and equally in each extended
fiber) can also be computed in polynomial time, if the dimension
is kept fix. It is an open question whether a similar short
rational function expression for the (right-hand sides of the)
atomic fibers exists and if it is polynomial time computable in
fixed dimension, that is, we ask for a short rational function
that encodes the generating functions
\[
\sum_{b\in F(A)} z^b\;\;\;\text{and}\;\;\;\sum_{b\in E(A)} z^b.
\]
At this point we would like to mention that one may of course use
Algorithm \ref{Algorithm to compute atomic fibers} also for the
problem of finding indecomposable extended fibers (and thus also
indecomposable fibers) among the elements in the family of
extended fibers $Q^I_{A,b}=\{z:Az=b,z\in\Z^n\}$ where $A$ is kept
fix and where $b$ is allowed to vary only on a sublattice
$\Lambda$ of $\Z^d$. One only needs to replace the input set
$F=\{Ae_1,\ldots,Ae_n,-Ae_1,\ldots,-Ae_n\}$ by a symmetric
generating set $F=\{\pm b_1,\ldots,\pm b_k\}$ of $\Lambda$ with
the additional property that $Q^I_{A,b}=\{z:Az=b,z\in\Z^n\}$ is
not empty for all $b\in F$.

{\bf Acknowledgment.} The author wishes to thank Jes\'us De Loera
for some useful suggestions on improving the presentation of this
paper.

\end{document}